\documentclass[11pt]{article}
\usepackage[latin1]{inputenc}
\usepackage[T1]{fontenc}
\usepackage{amsmath}
\usepackage{bbm}
\usepackage{amsthm}
\usepackage{amssymb,mathrsfs}
\usepackage{fourier}
\usepackage[english,french]{babel} 
\usepackage[all]{xy}
\usepackage{setspace}
\usepackage{color}
\usepackage{tabularx}

\newtheorem{theorem}{Theorem}[section]

\newtheorem{lemma}[theorem]{Lemma}
\newtheorem{proposition}[theorem]{Proposition}

\allowdisplaybreaks[1]

\newcommand\Z{{\mathbb Z}}

\newcommand\R{{\mathbb R}}
\newcommand\C{{\mathbb C}}

\begin{document}

\title{A short proof of the zero-two law for cosine functions}

\date{\relax}
\author{Jean Esterle}

\maketitle

{\bf {Abstract}}: Let $(C(t))_{t\in \R}$ be a cosine function in a unital  Banach algebra. We give a simple proof of the fact  that if lim sup$_{t\to 0}\Vert  C(t)-1_A\Vert<2,$ then  lim sup$_{t\to 0}\Vert  C(t)-1_A\Vert=0.$

{\it Keywords: Cosine function, scalar cosine function, commutative local Banach algebra.

}
\section{Introduction}

 Recall that a cosine function taking values in  a unital normed algebra $A$ with unit element $1_A$ is a family $C=(C(t))_{t\in \R}$ of elements of $A$  satisfying the so-called d'Alembert equation

\begin{equation}C(0)=1_A, C(s+t)+C(s-t)=2C(s)C(t)  \ \ (s\in \R,  t\in \R).\end{equation}

One can define in a similar way cosine sequences $(C(n))_{n\in \Z}.$ A cosine sequence depends only on the value of $C(1),$ since we have, for $n\ge 1,$

$$C(-n)=C(n)=T_n(C(1)),$$

where  $T_n(x)=\sum \limits _{k=0}^{[n/2]}C_n^kx^{n-2k}(x^2-1)^k$ is the $n^{th}$-Tchebyshev polynomial.

Strongly continuous operator valued cosine functions play an important role in the study of  abstract nonlinear second order differential equations, see for example \cite{tw}.
In a paper to appear in the Journal of Evolution Equations \cite{sz1}, Schwenninger and Zwart showed that if a strongly continuous cosine family $(C(t))_{t\in \R}$ of bounded operators on a Banach space $X$ satisfies lim sup$_{t\to 0}\Vert C(t) -I_X\Vert<2,$ then the generator $a$ of this cosine function is a bounded operator, so that lim sup$ \Vert C(t)-I_X\Vert=$lim sup$_{t\to 0}\Vert cos(ta)I_X-I_X\Vert=0,$ and they asked whether a similar zero-two law holds for general cosine functions $(C(t))_{t\in \R}.$ This question was answered positively by Chojnacki in \cite{ch1}. Using a sophisticated argument based on ultrapowers, Chonajcki deduced this zero-two law from the fact that if a cosine sequence $C(t)$ satisfies sup$_{t\in \R}\Vert C(t)-1_A\Vert <2,$ then $C(t)=1_A$ for $t \in \R.$ This second result, which was obtained independently by the author in \cite{e}, is proved by Chojnacki in \cite{ch1}
by adapting methods used by Bobrowski, Chojnacki and Gregosiewicz in \cite{bcg} to show that if a cosine sequence $(C(t))_{t\in \R}$ satisfies $sup_{t\in \R}\Vert C(t)-cos(at)1_A\Vert<{8\over 3\sqrt 3}$ for some $a\in \R,$ then $C(t)=cos(at)1_A$ for $t\in \R,$ a result also obtained independently by the author in \cite{e}, which improves previous results of \cite{bc}, \cite{ch2} and \cite{sz2}.

The purpose of this paper is to give a short direct proof of the zero-two law. The zero-two law for complex-valued cosine functions is a folklore result, which easily implies that if lim sup$_{t\to 0} \rho (C(t)-1_A)<2$ then lim sup ${t\to 0}\rho(C(t)-1_A)=0,$ where $\rho(x)$ denotes the spectral radius of an element $x$ of a Banach algebra, see section 2. Our proof of he zero-two law  is then based on the fact that if $\Vert C(t)-1_A\Vert \le 2,$ and if $\rho\left ( C \left ({t\over 2}\right )-1_A\right )<1,$ then we have

$$C\left ({t\over 2}\right )=\sqrt{1_A-{C(t)-1_A\over 2}},$$

where $\sqrt{1_A-u}$ is defined by the usual series for $\Vert u \Vert \le 1.$ It follows from this identity and from the fact that the coefficients of the Taylor series at the origin of the function $t \to 1-\sqrt{1-t}$ are positive  that in this situation we  have

$$\left \Vert C\left ({t\over 2}\right )-1_A\right \Vert \le 1 -\sqrt{ 1 - \left \Vert {C(t)-1_A\over 2}\right \Vert}, $$

and the zero-two law follows.

Notice that if we replace the constant $2$ by ${3\over 2}$ a "three line argument" due to Arendt \cite{a} shows that if lim sup$_{t\to 0} \Vert C(t)-1_A\Vert <{3\over 2}$ then 
lim sup$_{t\to 0} \Vert C(t)-1_A\Vert =0.$ The proof presented here has some analogy with Arendt's proof, and the difficulty to estimate $\Vert \left ( 1_A+C\left ({t\over 2}\right )\right )^{-1}\Vert$ is circumvented in the present paper by using the formula above.

The author wishes to thank W. Chonajcki for giving information about the reference  \cite{ch1}. He also thanks F. Schwenninger for valuable discussions about the content of the paper.

\section{The zero-two law for the spectral radius}

The zero-two law for scalar cosine functions pertains to folklore, but we could not find a reference in the litterature for the following certainly well-known lemma, which is a variant 
 of proposition 3.1 of \cite{e}.

\begin{lemma} Let $c=(c(t))_{t\in \R}$ be a complex-valued cosine function. Then $c$ statisfies one of the following conditions

(i)  $lim sup_{t\to 0}\vert c(t)-1\vert=+\infty,$

(ii)  $lim sup_{t\to 0}\vert c(t)-1\vert=2,$

(iii)  $lim sup_{t\to 0}\vert c(t)-1\vert=0.$
\end{lemma}

First assume that $M:= lim sup_{t\to 0}\vert c(t)\vert <+\infty,$ and denote by $S$ the set of all complex numbers $\alpha$ for which there exists a sequence $(t_m)_{m\ge 1}$ of positive real numbers such that $lim_{m\to +\infty}t_m=0$ and $lim_{m\to +\infty}c(t_m)=\alpha.$ Then $\vert \alpha \vert \le M$ for every $\alpha \in S.$ Notice that if $\alpha \in S,$ and if a sequence $(t_m)_{m\in Z}$ satisfies the above conditions with respect to $\alpha,$ then $T_n(\alpha)= lim_{m\to +\infty}T_n(C(t_m))=lim_{m\to +\infty}C(nt_m),$ and so $T_n(\alpha)\in S,$ and  $\vert T_n(\alpha)\vert \le M$ for $n\ge 1.$ Now write $\alpha =cos(z)=\sum \limits_{k=0}^{+\infty}{z^{2k}\over (2k)!},$ and set $u=Re(z), v=Im(z).$ We have, for $n\ge 1,$

$$T_n(\alpha)=cos(nz)={e^{inu}e^{-nv}+e^{-inu}e^{nv}\over 2}.$$

Since $sup_{n\ge 1}\vert T_n(\alpha)\vert \le M,$ we have $v=0,$ $S\subset [-1,1],$ and $lim sup_{t\to 0}\vert c(t)-1\vert \le 2.$

Now assume that $S\neq \{1\},$ and let $\alpha \in S\setminus \{1\}.$ We have $\alpha =cos(u)$ for some $u \in \R.$ We see as above that $cos(nu)\in S$ for every $n \ge 1.$ If $u/\pi$
is irrationnal, then the set $\{e^{inu}\}_n\ge 1$ is dense in the unit circle $\mathbb T,$ and so $S=[-1,1]$ since $S$ is closed, and in this situation $lim sup_{t\to 0}\vert c(t)-1\vert=2.$

Now assume that $u/ \pi$ is rational, and let $s \ge 2$ be the smallest positive integer such that $e^{ius}=1.$ Then $e^{2i\pi \over s}=e^{ipu}$ sor some positive integer $p,$ and so $cos\left ({2\pi\over s}\right ) \in S.$  Let $(t_m)_{m\ge 1}$ be a sequence of positive reals such that $lim_{m\to +\infty}t_m=0$ and $lim_{m\to +\infty}c(t_m)=cos\left ({2\pi\over s}\right ),$ let $q\ge 2,$ and let $\beta$ be a limit point of the sequence $c\left ({t_m\over s^{q-1}}\right )_{n\ge 1}.$ There exists $y \in \R$ such that $cos(y)=\beta$ and such that $s^{q-1}y= {2\pi\over s} +{2k\pi},$ with $k \in \Z.$ Then $y=(1+ks){2\pi\over s^q}.$ Since $gcd(1+ks,s^q)=1,$  there exists a positive integer $r$ such that $ry - {2\pi\over s^q} \in 2\pi\Z,$ so that  $cos\left ({2\pi\over s^q}\right )\in S.$ This implies that $cos\left ( {2p\pi\over s^q}\right )\in S$ for $p\ge 1, q \ge 1,$ and $S$ is dense in $[-1,1].$ Since $S$ is closed, we obtain again $S=[-1,1],$ which implies that $lim sup_{t\to 0}\vert c(t)-1\vert=2.$ So if neither (i) nor (ii) holds, we have $S=\{1\},$ which implies (iii).  $\square$

Notice that if a cosine function $C=(C(t))_{t\in \R}$ in a Banach algebra $A$ satisfies $sup_{\vert t \vert \le \eta}\Vert C(t)\Vert \le M <+\infty$ for some $\eta >0,$ then $sup_{\vert t \vert \le L}\Vert C(t)\Vert <+\infty$ for every $L>0,$ since $sup _{\vert t \vert \le n\eta}\Vert C(t)\Vert \le sup_{\Vert y \Vert \le M}\Vert T_n(y)\Vert$ for every $n \ge 1,$ where $T_n$ denotes the $n$-th Tchebyshev polynomial. In particular if a complex-valued cosine function $c=(c(t))_{t\in \R}$ satisfies (iii), then the identity

$$(1-c(s-t))(1-c(s+t))=(c(s)-c(t))^2$$

shows as is well-known that the cosine function $c$ is continuous on $\R,$ which implies that $c(t)=cos(ta)$ for some $a \in \C.$

 If $A$ is commutative and unital, we will denote $1_A$ the unit element of $A,$ and we will denote by $\widehat A$ the space of all characters on $A,$ equipped with the Gelfand topology, i.e. the compact topology induced by the weak$*$ topology on the unit ball of the dual space of $A.$

\begin{proposition} Let $C=(C(t))$ be a cosine function in a unital Banach algebra $A.$ If $lim sup_{t\to 0}\rho (C(t)-1_A)<2,$ then $lim sup_{t\to 0}\rho (C(t)-1_A)=0.$

\end{proposition}

Proof: We may assume that unital Banach algebra $A$ is genarated by $(C(t))_{t\in \R}.$ Let $\chi \in \widehat A.$ Then the cosine complex-valued function $(\chi(C(t)))_{t\in \R}$ satisfies condition (iii) of the lemma, and so there exists $a_{\chi}\in \C$ such that we have

$$\chi(C(t))=cos(ta_{\chi}) \ \ (t \in \R).$$ 

Set $u_{\chi}=Re(a_{\chi}), v_{\chi}=Im(a_{\chi}).$ We have

$$\rho(C(t)-1)\ge \vert 1 -cos(tu_{\chi})ch(tv_{\chi})\vert.$$

If the family $(u_{\chi})_{\chi \in \widehat A}$ were unbounded, there would exist a sequence $(t_n)_{n\ge 1}$ of real numbers converging to zero and a sequence $(\chi_n)_{n\ge 1}$ of characters of $A$ such that $cos(t_nu_{\chi_n})=-1,$ and we would have $\rho(C(t_n)-1)\ge 2$ for $n \ge 1.$ So the family $(u_{\chi})_{\chi \in \widehat A}$ is bounded. If the family $(v_{\chi})_{\chi \in \widehat A}$ were unbounded, there would exist a sequence $(t'_n)_{n\ge 1}$ of real numbers converging to zero and a sequence 
$(\chi'_n)_{n\ge 1}$ of characters of $A$ such that $lim _{n\to +\infty}ch(t'_nv_{\chi'_n})=+\infty.$ But this would imply that $lim sup_{t\to 0}\rho (C(t))=+\infty.$ 
Hence the family $(a_{\chi})_{\chi \in \widehat A}$ is bounded, and we have

$$lim sup_{t\to 0}\rho(C(t)-1_A)=lim_{t\to 0}sup_{\chi \in \widehat A}\vert cos(ta_{\chi})-1\vert =0.$$ $\square$

\section{The zero-two law for cosine functions}

Set $\alpha_n={1\over n!}{1\over 2}\left ({1\over 2}-1\right )\dots \left ({1\over 2}-n+1\right )$ for $n\ge 1,$ with the convention $\alpha_{0}=0,$ and for $\vert z \vert <1,$ set

$$\sqrt{1-z}=\sum_{n=0}^{+\infty}(-1)^n\alpha_nz^n,$$

so that $\left (\sqrt{1-z}\right )^2=1-z,$ and $\sqrt{1-t}$ is the  positive square root of $1-t$ for $t\in (-1,1).$ Also $Re\left (\sqrt {1-z}\right )>0$ for $\vert z \vert <1.$

 We have, for $t \in [0,1),$

$$\sum_{n=1}^{+\infty}(-1)^{n-1}\alpha_nt^n=1-\sqrt {1-t}.$$

Since $(-1)^{n-1}\alpha_n\ge 0$ for $n\ge 1,$ the series $\sum_{n=1}^{+\infty}\vert \alpha_n\vert =\sum_{n=1}^{+\infty}(-1)^{n-1}\alpha_n$ is convergent, and we have

$$\sum_{n=1}^{+\infty}\vert \alpha_n\vert t^n=1-\sqrt {1-t} \ \ (0\le t \le 1).$$

Now let $A$ be a commutative unital Banach algebra, and let $x \in A$ such that $\Vert x \Vert\le 1.$ Set

$$\sqrt{1_A-x}=\sum_{n=0}^{+\infty}(-1)^n\alpha_nx^n.$$

Then $\left (\sqrt{1_A-x}\right )^2=1_A-x,$  and we have

\begin{equation}\left \Vert 1_A-\sqrt{1_A-x}\right \Vert =\left \Vert \sum_{n=1}^{+\infty}(-1)^n\alpha_n x^n\right \Vert \le \sum_{n=1}^{+\infty}\vert \alpha_n\vert \Vert  x\Vert ^n=1-\sqrt{1-\Vert x
\Vert}
\end{equation}

Notice also that if $A$ is commutative, then we have

\begin{equation} Re\left ( \chi \left (\sqrt{1_A -x}\right ) \right )=Re \left ( \sqrt{1-\chi(x)}\right ) \ge 0 \  \ (\chi \in \widehat A).\end{equation}

We obtain the following formula

\begin{lemma} Let $(C(t))_{t\in \R}$ be a cosine function in a unital Banach algebra $A.$ Assume that $\Vert C(t)-1_A\Vert \le 2$ and that $\rho(\left (C\left ({t\over 2}\right )-1\right )<1.$ Then we have

$$C\left ({t\over 2}\right )=\sqrt { 1_A-{1_A-C(t)\over 2}}.$$

\end{lemma}

Proof: The abstract version of the formula $sin^2\left ({u\over 2}\right )={1-cos(u)\over 2}$ gives

$$ 1_A -C\left ({t\over 2}\right )^2= {1_A-C(t)\over 2}, C\left ({t\over 2}\right )^2=1_A-{1_A-C(t)\over 2}=\left (\sqrt { 1_A-{1_A-C(t)\over 2}}\right )^2,$$

$$\left (C\left ({t\over 2}\right )-\sqrt { 1_A-{1_A-C(t)\over 2}}\right )\left (C\left ({t\over 2}\right )+\sqrt { 1_A-{1_A-C(t)\over 2}}\right )=0.$$

We may assume that $A$ is commutative. Let $\chi \in \widehat A.$ Since $\rho\left (C\left ({t\over 2}\right )-1_A\right )<1,$ we have $Re\left ( \chi  \left (C\left ({t\over 2}\right )\right )\right )>0.$ Since $Re\left (\chi \left ( \sqrt { 1_A-{1_A-C(t)\over 2}}\right )\right )\ge 0,$ $C\left ({t\over 2}\right )+\sqrt { 1_A-{1_A-C(t)\over 2}}$ is invertible in $A,$ and $C\left ({t\over 2}\right )-\sqrt { 1_A-{1_A-C(t)\over 2}}=0.$ $\square$

\begin{theorem} Let $(C(t))_{t\in \R}$ be a cosine sequence in a Banach algebra. If $lim sup_{t\to 0}\Vert C(t)-1_A\Vert <2,$ then $lim sup_{t\to 0}\Vert C(t)-1_A\Vert =0.$

\end{theorem}

Proof: It follows from proposition 2.2 and lemma 3.1 that there exists $\eta >0$ such that we have, for $\vert t \vert \le \eta,$

$$ \Vert C(t)-1_A\Vert <2, C\left ({t\over 2}\right )=\sqrt { 1_A-{1_A-C(t)\over 2}}.$$

Using (1), we obtain, for $\vert t \vert \le \eta,$

$$\left \Vert C\left ({t\over 2}\right )-1_A\right \Vert \le 1 -\sqrt{ 1 - \left \Vert {C(t)-1_A\over 2}\right \Vert}.$$

Set $l=lim sup_{t\to 0}\Vert C(t)-1_A\Vert.$ We obtain

$$l \le 1 -\sqrt{1-{l\over 2}} \le 1,$$

and so $l=0.$ $\square$

Notice that the proof above gives a little bit more than the zero-two law: if $\Vert 1-C(t)\Vert \le 2$ and $\rho\left (1-C\left ({t\over 2}\right )\right )<1$ for $\vert t \vert \le \eta,$ then we have, for $n\ge 1,$

$$sup_{\vert t \vert \le 2^{-n}\eta}\Vert C(t)-1_A\Vert \le u_n,$$

where the sequence $u_n$ satisfies $u_0=2, u_{n+1}= 1-\sqrt{1-{u_n\over 2}}$ for $n\ge 1,$ and $lim_{n\to +\infty}u_n=0,$ which gives an explicit control on the convergence to $0$ of $\Vert C(t)-1_A\Vert $ as $t\to 0.$

Notice also that the fact that the coefficients of the Taylor expansion at the origin of the function $t \to 1-\sqrt{1-t}$ are positive was used in \cite{e0} to show that $\Vert x^2-x\Vert \ge 1/4$ for every quasinilpotent element $x$ of a Banach algebra such that $\vert x\vert \ge 1/2.$ Similar argument were used in \cite{em} to show that if a semigroup $(T(t))_{t>0}$ in a Banach algebra $A$ satisfies lim sup$_{t\to 0^+}\Vert T(t) -T((n+1)t)\Vert < {n\over (n+1)^{1+{1\over n}}}$ for some $n \ge 1$, then there exists an idempotent $J$ of $A$ such that lim$_{t\to 0}\Vert T(t)-J\Vert=0,$ so that  lim sup$_{t\to 0^+}\Vert T(t) -T((n+1)t)\Vert =0.$

{\it Jean Esterle

IMB, UMR 5251

Universit\'e de Bordeaux

351, cours de la Lib\'eration

33405 -Talence (France)

esterle@math.u-bordeaux1.fr}

\end{document}